\documentclass[12pt,a4paper,reqno]{article}%
\usepackage{amssymb}
\usepackage{amsmath}
\usepackage{amsfonts}
\usepackage{graphicx}
\usepackage{mathrsfs}
\usepackage[all]{xy}
\usepackage{hyperref}%
\providecommand{\U}[1]{\protect\rule{.1in}{.1in}}
\SelectTips{cm}{}
\newtheorem{theorem}{Theorem}

\newtheorem{corollary}[theorem]{Corollary}

\begin{document}

\title{\textbf{On Higher Derivatives as Constraints in Field Theory: a Geometric
Perspective}}
\author{\textsc{{L.~Vitagliano}\thanks{\textbf{e}-\textit{mail}:
\texttt{lvitagliano@unisa.it}}}\\{\small {DMI, Universit\`{a} degli Studi di Salerno, Via Ponte don Melillo,
84084 Fisciano (SA), Italy}}\\{\small INFN, CG di Salerno - Sezione di Napoli, via Cintia, 80126 Naples,
Italy}\\{\small Levi-Civita Institute, via Colacurcio 54, 83050 Santo Stefano del Sole
(AV), Italy}}
\maketitle

\begin{abstract}
We formalize geometrically the idea that the (de Donder) Hamiltonian
formulation of a higher derivative Lagrangian field theory can be constructed
understanding the latter as a first derivative theory subjected to constraints.

\end{abstract}

\emph{Keywords}: Higher Derivative Field Theory, Fiber Bundles, Jet Bundles,
Lagrangian and Hamiltonian Formalisms, Constraints.

\quad

\emph{2000 MSC}: 53B50, 53C80, 70S05.

\section{Introduction}

Let $\pi:E\longrightarrow M$ be a fiber bundle, $\pi_{l}:J^{l}\longrightarrow
M$ its $l$-th jet bundle, $l=0,1,2,\ldots$, and $\mathscr{L}\in\Lambda
^{n}(J^{k})$ a basic $n$-form on $J^{k}$, $n=\dim M$, $k>1$. $\mathscr{L}$ may
be interpreted as a Lagrangian density defining the $k$th derivative action
functional $\boldsymbol{S}:s\longmapsto\int_{M}(j_{k}s)^{\ast}\mathscr{L}$ on
sections $s$ of $\pi$. The associated calculus of variations, and, in
particular, the Euler-Lagrange equations, have a nice geometric (and
homological) formulation in terms of the so called $\mathscr{C}$-spectral
sequence \cite{v84}. The Hamiltonian counterpart of the theory is very well
established in the case $k=1$. In particular, there are universally accepted
field theoretic, geometric definitions of the Legendre transform and the
Hamilton equations (see, for instance, \cite{r09} for a recent review). On the
other hand, only recently a geometric formulation of the (Hamilton-like,
higher derivative) de Donder field theory \cite{d35} has been proposed by the
author which is natural, i.e. it is independent on any structure other than
the action functional itself \cite{v10}. Such a formulation is based on a
generalization to higher derivative Lagrangian field theory of the mixed
Lagrangian-Hamiltonian formalism by Skinner and Rusk \cite{s83,sr83,sr83b}. In
such a theory the Legendre transform is not defined a priori but it is rather
a consequence of the field equations.

Aldaya and de Azcarraga have suggested that higher derivative Hamiltonian
field theory can be introduced understanding higher derivative Lagrangian
field theory as a first order theory with (vakonomic) constraints \cite{aa82}.
However, they work in local coordinates and not all their conclusions have an
intrinsic, geometric meaning. The aim of this short communication is to show
that the idea by Aldaya and de Azcarraga can be given a precise, and natural,
geometric formulation. In particular, momenta in higher derivative field
theory can be mathematically understood as Lagrange multipliers in an
equivalent first derivative theory subjected to (vakonomic) constraints.

\section{The Constraint Bundle}

We assume that the reader is familiar with Lagrangian and Hamiltonian
formalisms on fiber bundles \cite{r09}. We refer to \cite{v10} and \cite{v09}
for notations, conventions, and the basic differential geometric constructions
we will use in the following.

Let $\mathscr{L}$ be as in the introduction. It is well known that $J^{k}$ is
naturally embedded in $J^{1}\pi_{k-1}$, the first jet bundle of $\pi_{k-1}$.
Denote by $\ldots,x^{i},\ldots$ coordinates on $M$, by $\ldots,u_{I}^{\alpha
},\ldots$ jet coordinates on $J^{k}$, ($I=i_{1}\cdots i_{r}$ being a
multi-index denoting multiple differentiation of the field variables
$\ldots,u^{\alpha},\ldots$, $i_{1},\ldots i_{r}=1,\ldots,n$, $|I|{}:=r\leq k$)
and by $\ldots,u_{J}^{\alpha}{}_{.i},\ldots$ jet coordinates on $J^{1}%
\pi_{k-1}$, $|J|{}\leq k-1$. The embedding $e:J^{k}\hookrightarrow J^{1}%
\pi_{k-1}$ reads locally
\[
e^{\ast}(u_{J}^{\alpha}{}_{.i})=u_{Ji}^{\alpha},\quad|J|{}\leq k-1.
\]
In particular $J^{k}\longrightarrow J^{k-1}$ is an affine subbundle of
$J^{1}\pi_{k-1}\longrightarrow J^{k-1}$. $\mathscr{L}$ can be understood as a
first derivative Lagrangian density, say $\mathscr{L}^{\prime}$, on the
constraint subbundle $J^{k}$ of $J^{1}\pi_{k-1}$. Sections $\sigma$ of
$\pi_{k-1}$ satisfying the constraint, i.e. whose first jet prolongation
$j_{1}\sigma$ takes values in $J^{k}\subset J^{1}\pi_{k-1}$, are precisely
those of the form $\sigma=j_{k-1}s$ for some section $s$ of $\pi$. In other
words, considering $J^{k}$ as a constraint subbundle of $J^{1}\pi_{k-1}$ is
the same as introducing new variables corresponding to derivatives of the
fields and then impose the obvious differential relations among them.
Therefore, the variational problem defined by $\mathscr{L}^{\prime}$ is
equivalent to the original one, and, in principle, we can apply the Lagrange
multiplier method to find solutions. To do this, we should, first of all, 1)
choose an extension of $\mathscr{L}^{\prime}$ to the whole $J^{1}\pi_{k-1}$
and 2) present $J^{k}\subset J^{1}\pi_{k-1}$ as the zero locus of a
(sufficiently regular) morphism of the bundle $J^{1}\pi_{k-1}\longrightarrow
J^{k-1}$, with values in a vector bundle $V\longrightarrow J^{k-1}$
\cite{ggr06}. Since neither 1) nor 2) can be done in a natural way, we prefer
to change a bit our strategy.

Instead of $J^{1}\pi_{k-1}$, consider $J^{1}\pi_{k}$, the first jet bundle of
$\pi_{k}$. There is a natural projection $p:J^{1}\pi_{k}\longrightarrow J^{k}%
$. Moreover, we can draw a diagram
\[%
\begin{array}
[c]{c}%
\xymatrix{   X_{k}\ \ar@<-0.5ex>@{^(->}
[r]  \ar[d]  &         J^1 \pi_{k} \ar[d] \\
J^{k}\ \ar@<-0.5ex>@{^(->}[r]           &             J^1 \pi_{k-1}}%
\end{array}
,
\]
where $X_{k}:=p^{-1}(J^{k})$. We understand $X_{k}$ as a contraint subbundle
in $J^{1}\pi_{k}$. Notice that $\mathscr{L}$ is naturally extended to
$J^{1}\pi_{k}$ (and, in particular, $X_{k}$) as $p^{\ast}(\mathscr{L})$.
Moreover, $X_{k}\subset J^{1}\pi_{k}$ can be presented as the zero locus of a
morphism $\psi:J^{1}\pi_{k}\longrightarrow V$ of the bundle $J^{1}\pi
_{k}\longrightarrow J^{k}$, with values in a vector bundle $V\longrightarrow
J^{k}$, as follows. Let $\theta_{0}\in J^{k}$ and $\theta\in J^{1}\pi_{k}$ be
a point over it, i.e., the projection $J^{1}\pi_{k}\longrightarrow J^{k}$
sends $\theta$ to $\theta_{0}$. $\theta$ can be understood as an
$n$-dimensional subspace $L(\theta)$ in $T_{\theta_{0}}J^{k}$ transversal to
the fiber $F$ of $\pi_{k}$ through $\theta_{0}$ (see, for instance,
\cite{b...99}), or, which is the same, as a linear map $\Pi(\theta
):T_{\theta_{0}}J^{k}\longrightarrow T_{\theta_{0}}J^{k}$, with the following
two properties: 1) $\Pi(\theta)$ is a projector, i.e., $\Pi(\theta)\circ
\Pi(\theta)=\Pi(\theta)$, 2) $\ker\Pi(\theta)=T_{\theta_{0}}F$. Then
$L(\theta)=\operatorname{im}\Pi(\theta)$. If $\theta$ has jet coordinates
$\ldots,u_{I}^{\alpha}{}_{.i},\ldots$, $|I|{}\leq k$, then
\[
\Pi(\theta)=\left(  \frac{\partial}{\partial x^{i}}+\sum\nolimits_{|I|{}\leq
k}u_{I}^{\alpha}{}_{.i}\frac{\partial}{\partial u_{I}^{\alpha}}\right)
\otimes dx^{i}.
\]

Now, there is a canonical geometric structure on $J^{k}$, the so called Cartan
distribution \cite{b...99}. The Cartan plane $\mathscr{C}(\theta_{0})\subset
T_{\theta_{0}}J^{k}$ at $\theta_{0}$ can be described as the kernel of a
canonical linear map $U(\theta_{0}):T_{\theta_{0}}J^{k}\longrightarrow
T_{\bar{\theta}_{0}}J^{k-1}$, $\bar{\theta}_{0}\in J^{k-1}$ being the image of
$\theta_{0}$ under the projection $J^{k}\longrightarrow J^{k-1}$ \cite{klv86}.
In local coordinates $U(\theta_{0})$ is given by
\[
U(\theta_{0})=\sum\nolimits_{|I|{}\leq k-1}\frac{\partial}{\partial
u_{I}^{\alpha}}\otimes(du_{I}^{\alpha}-u_{Ii}^{\alpha}dx^{i}).
\]
We can also compose $\Pi(\theta)$ and $U(\theta_{0})$, to check wether
$L(\theta)\subset\mathscr{C}(\theta_{0})$. In local coordinates
\begin{align}
U(\theta_{0})\circ\Pi(\theta)  &  =\sum\nolimits_{|I|{}\leq k-1}%
(du_{I}^{\alpha}-u_{Ii}^{\alpha}dx^{i})\left(  \frac{\partial}{\partial x^{j}%
}+\sum\nolimits_{|J|{}\leq k}u_{J}^{\beta}{}_{.j}\frac{\partial}{\partial
u_{J}^{\beta}}\right)  \frac{\partial}{\partial u_{I}^{\alpha}}\otimes
dx^{j}\nonumber\\
&  =\sum\nolimits_{|I|{}\leq k-1}(u_{I}^{\alpha}{}_{.i}-u_{Ii}^{\alpha}%
)\frac{\partial}{\partial u_{I}^{\alpha}}\otimes dx^{i}. \label{Eq1}%
\end{align}
We conclude that $L(\theta)\subset\mathscr{C}(\theta_{0})$ iff $\ldots
,u_{I}^{\alpha}{}_{.i}=u_{Ii}^{\alpha},\ldots$, $|I|{}\leq k-1$, i.e.,
$\theta\in X_{k}$.

In view of its coordinate expression, $U(\theta_{0})\circ\Pi(\theta)$ can be
understood as an element in $V_{\bar{\theta}_{0}}J^{k-1}\otimes T_{x}^{\ast}%
M$, where $V_{\bar{\theta}_{0}}J^{k-1}=\ker d_{\bar{\theta}_{0}}\pi
_{k-1}\subset T_{\bar{\theta}_{0}}J^{k-1}$ is the $\pi_{k-1}$-vertical tangent
space to $J^{k-1}$ at the point $\bar{\theta}_{0}$ and $x=\pi_{k}(\theta)\in
M$. Therefore, the map $\theta\longmapsto U(\theta_{0})\circ\Pi(\theta)$ can
be understood as an affine morphism $\psi:J^{1}\pi_{k}\longrightarrow V$ of
the bundle $J^{1}\pi_{k}\longrightarrow J^{k}$, with values in the (pull-back)
vector bundle
\[
V:=VJ^{k-1}\otimes_{M}T^{\ast}M\times_{J^{k-1}}J^{k}\longrightarrow J^{k}%
\]
whose fiber over $\theta_{0}$ is $V_{\bar{\theta}_{0}}J^{k-1}\otimes
T_{x}^{\ast}M$. Formula (\ref{Eq1}) then shows that $\theta\in X_{k}$ iff
$\psi(\theta)=0$. Formula (\ref{Eq1}) also shows that $\psi$ has fiber-wise
maximal rank at the points of $X_{k}$, and in this sense, will be referred to
as a \emph{regular morphism }\cite{ggr06}. We have thus proved the following

\begin{theorem}
\label{Theor1}$X_{k}\subset J^{1}\pi_{k}$ is the zero locus of a canonical
regular morphism of the affine bundle $J^{1}\pi_{k}\longrightarrow J^{k}$ with
values in a canonical vector bundle $V\longrightarrow J^{k}$.
\end{theorem}

Notice that $VJ^{k-1}\otimes_{M}T^{\ast}M\longrightarrow J^{k-1}$ is the model
vector bundle for the affine bundle $J^{1}\pi_{k-1}\longrightarrow J^{k-1}$.

\begin{corollary}
A smooth function $F\in C^{\infty}(J^{1}\pi_{k})$ vanishes on the constraint
subbundle $X_{k}$ iff there exists a morphism $\lambda:J^{1}\pi_{k}%
\longrightarrow V^{\ast}$, with values in the dual bundle, such that
$\langle\lambda,\psi\rangle=0$.
\end{corollary}

The above corollary shows that variables in the fiber of $V^{\ast
}\longrightarrow J^{k}$ basically play the role of Lagrange multipliers (see
below for details).

\section{Higher Derivatives as Constraints}

Consider the first derivative action functional $\boldsymbol{S}^{\prime
}:\sigma\longmapsto\int_{M}\sigma^{\ast}\mathscr{L}$ on sections of $\pi_{k}$
constrained by $X_{k}$, i.e., we restrict $\boldsymbol{S}^{\prime}$ to those
sections $\sigma$ such that $\operatorname{im}j_{1}\sigma\subset X_{k}$
(notice that, without the constraints, $\boldsymbol{S}^{\prime}$ would
actually be a zeroth derivative action functional and, therefore, a very
trivial one). The variational problem defined in this way is equivalent to the
original one. In fact, similarly as above, sections $\sigma$ of $\pi_{k}$ such
that $\operatorname{im}j_{1}\sigma\subset X_{k}$ are precisely those of the
form $\sigma=j_{k}s$ for some section $s$ of $\pi$. In view of Theorem
\ref{Theor1}, we can use the method of Lagrange multipliers to find extremals.
In the present case, the method consists in searching for extremals of a new,
unconstrained, first derivative, action functional $\boldsymbol{S}_{1}%
:\Sigma\longmapsto\int_{M}(j_{1}\Sigma)^{\ast}\mathscr{L}_{1}$ on an augmented
space of sections $\Sigma$. More precisely, $\Sigma$ is a section of the
bundle $V^{\dag}:=V^{\ast}\otimes_{M}\Lambda^{n}T^{\ast}M\longrightarrow M$,
which, in the following, we denote by $q$. Notice that, by construction,
points of $V$ and points of $V^{\dag}$ over the same point $\theta_{0}$ of
$J^{k}$ can be paired to give a top form over $M$ at $\pi_{k}(\theta_{0})$. We
denote by $\langle\cdot,\cdot\rangle$ such pairing. Since $\operatorname{Hom}%
(\Lambda^{1}(M),\Lambda^{n}(M))\simeq\Lambda^{n-1}(M)$ we have
\[
V^{\dag}\simeq V^{\ast}J^{k-1}\otimes_{M}\Lambda^{n-1}T^{\ast}M\times
_{J^{k-1}}J^{k}%
\]
and it identifies naturally with $J^{\dag}\pi_{k-1}\times_{J^{k-1}}J^{k}$,
$J^{\dag}\pi_{k-1}$ being the reduced multi-momentum bundle of $\pi_{k-1}$
\cite{r09} (see also \cite{v10}).

The Lagrangian density $\mathscr{L}_{1}$ is defined by
\[
(j_{1}\Sigma)^{\ast}\mathscr{L}_{1}=\sigma^{\ast}\mathscr{L}+\langle
\Sigma,\psi\circ j_{1}\sigma\rangle\in\Lambda^{n}(M),
\]
where $\Sigma$ is a section of $q:V^{\dag}\longrightarrow M$, and $\sigma$ is
the section of $\pi_{k}$ given by projecting $\Sigma$ onto $J^{k}$. Describe
$\mathscr{L}_{1}$ locally. To this aim, let $\mathscr{L}$ be locally given by
\[
\mathscr{L}=L[x,u]d^{n}x,
\]
where $L[x,u]:=L(\ldots,x^{i},\ldots,u_{I}^{\alpha},\ldots)$, $|I|{}\leq k$,
is a local function on $J^{k}$ and $d^{n}x:=dx^{1}\wedge\cdots\wedge dx^{n}$.
Moreover, let $\ldots,p_{\alpha}^{J}{}^{.j},\ldots$be standard, dual
coordinates on $J^{\dag}\pi_{k-1}$ corresponding to jet coordinates
$\ldots,u_{J}^{\alpha},\ldots$ on $J^{k-1}$, $|J|{}\leq k-1$. It is easy to
see that, locally, $\mathscr{L}_{1}=L_{1}[x,u,p,u.,p.]d^{n}x$ where
\begin{align}
L_{1}[x,u,p,u.,p.]:  &  =L_{1}(\ldots,x^{i},\ldots,u_{I}^{\alpha}%
,\ldots,p_{\alpha}^{J}{}^{.j},\ldots,u_{I}^{\alpha}{}_{.i},\ldots,p_{\alpha
}^{J}{}^{.j}{}_{.i},\ldots)\nonumber\\
&  =L[x,u]+\sum\nolimits_{|I|{}\leq k-1}p_{\alpha}^{I}{}^{.i}(u_{I}^{\alpha}%
{}_{.i}-u_{Ii}^{\alpha}), \label{Eq2}%
\end{align}
$\ldots,p_{\alpha}^{J}{}^{.j}{}_{.i},\ldots$ being jet coordinates
corresponding to coordinates $\ldots,p_{\alpha}^{J}{}^{.j},\ldots$ on
$V^{\dag}$. Formula (\ref{Eq2}) shows that the $p_{\alpha}^{J}{}^{.j}$'s,
i.e., variables in the fiber of $V^{\dag}\longrightarrow J^{k}$ play the role
of Lagrange multipliers.

Now, consider the Euler-Lagrange-Hamilton equations \cite{v10} determined by
$\boldsymbol{S}$. They are the higher derivative, field theoretic analogue of
the equations of motions of a Lagrangian mechanical system proposed by Skinner
and Rusk in \cite{s83,sr83}. Recall that the Euler-Lagrange-Hamilton equations
are imposed precisely on sections of $V^{\dag}\longrightarrow M$ and are of
the PD-Hamilton type (see \cite{v09} for the definition and main properties of
PD-Hamiltonian systems and their PD-Hamilton equations). Moreover, the
PD-Hamiltonian system determining them is an exact form. Therefore, the
Euler-Lagrange-Hamilton equations are the Euler-Lagrange equations determined
by a suitable Lagrangian density. The latter coincides with $\mathscr{L}_{1}$
up to total divergences. Indeed, the Euler-Lagrange equations determined by
$\boldsymbol{S}_{1}$ locally read
\[
\left(
\begin{array}
[c]{c}%
\dfrac{\delta}{\delta u_{I}^{\alpha}}\\
\dfrac{\delta}{\delta p_{\alpha}^{J}{}^{.j}}%
\end{array}
\right)  L_{1}=0.
\]
Now,
\[
\left(
\begin{array}
[c]{c}%
\dfrac{\delta}{\delta u_{I}^{\alpha}}\\
\dfrac{\delta}{\delta p_{\alpha}^{J}{}^{.j}}%
\end{array}
\right)  L_{1}=\left(
\begin{array}
[c]{c}%
\dfrac{\partial}{\partial u_{I}^{\alpha}}-\dfrac{d}{dx^{i}}\dfrac{\partial
}{\partial u_{I}^{\alpha}{}_{.i}}\\
\dfrac{\partial}{\partial p_{\alpha}^{J}{}^{.j}}-\dfrac{d}{dx^{i}}%
\dfrac{\partial}{\partial p_{\alpha}^{J}{}^{.j}{}_{.i}}%
\end{array}
\right)  L_{1}=\left(
\begin{array}
[c]{c}%
\dfrac{\partial L}{\partial u_{I}^{\alpha}}-\delta_{Jj}^{I}p_{\alpha}^{J}%
{}^{.j}-p_{\alpha}^{I}{}^{.i}{}_{.i}\\
u_{J}^{\alpha}{}_{.j}-u_{Jj}^{\alpha}%
\end{array}
\right)
\]
which is the left hand side of Euler-Lagrange-Hamilton equations determined by
$\boldsymbol{S}$. Summarizing, we have proved the following

\begin{theorem}
The Euler-Lagrange equations determined by $\boldsymbol{S}_{1}$ coincide with
the Euler-Lagrange-Hamilton equations determined by $\boldsymbol{S}$.
\end{theorem}

Recall that the Euler-Lagrange-Hamilton equations cover the Euler-Lagrange
equations in the sense that solutions of the former are surjectively mapped to
solutions of the latter by projection onto $E$ \cite{v10}. We then duly
recover the Lagrange multiplier theorem in the present case (see, for
instance, \cite{ggr06}, see also \cite{mps01}).

\section{The Hamiltonian Sector}

Let us now have a look at the Hamiltonian counterpart of the field theory
defined by $\boldsymbol{S}_{1}$. Let $J^{\dag}q$ be the reduced multimomentum
bundle of $q:V^{\dag}\longrightarrow M$ and $\ldots,P_{\alpha}^{I}{}%
^{.i},\ldots,Q_{J}^{\alpha}{}{}_{.j}{}^{.i},\ldots$ be dual coordinates on it
corresponding to coordinates $\ldots,u_{I}^{\alpha},\ldots,p_{\alpha}^{J}%
{}^{.j},\ldots$ on $V^{\dag}$, respectively, $|I|{}\leq k$, $|J|{}\leq k-1$.
The Legendre transform $F\mathscr{L}_{1}:J^{1}q\longrightarrow J^{\dag}q$ is
the fiber-derivative of $\mathscr{L}_{1}$ \cite{r09}. Clearly,
$F\mathscr{L}_{1}$ is actually independent of $\mathscr{L}$. Locally,
\begin{align*}
F\mathscr{L}_{1}^{\ast}(P_{\alpha}^{I.i})  &  =\left\{
\begin{array}
[c]{ll}%
p_{\alpha}^{I.i} & \text{if }|I|{}\leq k-1\\
0 & \text{if }|I|{}=k
\end{array}
\right.  ,\\
F\mathscr{L}_{1}^{\ast}(Q_{J.j}^{\alpha}{}^{.i})  &  =0.
\end{align*}
This shows that $\operatorname{im}F\mathscr{L}_{1}\simeq V^{\dag}$ and that,
if we understand this isomorphism, $F\mathscr{L}_{1}:J^{1}q\longrightarrow
V^{\dag}$ is nothing but the canonical projection. In particular,
$F\mathscr{L}_{1}:J^{1}q\longrightarrow\operatorname{im}F\mathscr{L}_{1}$ is a
surjective submersion with connected fibers and, therefore, $\mathscr{L}_{1}$
induces on $\operatorname{im}F\mathscr{L}_{1}\simeq V^{\dag}$ a unique
PD-Hamiltonian system $\omega$ such that $F\mathscr{L}_{1}^{\ast}%
(\omega)=d\Theta_{\mathscr{L}_{1}}$, $\Theta_{\mathscr{L}_{1}}$ being the
Poincar\'{e}-Cartan $n$-form determined by $\mathscr{L}_{1}$ on $J^{1}q$
\cite{r09}. A direct computation shows that $\omega$ is locally given by
\[
\omega=\sum\nolimits_{|I|{}\leq k-1}dp_{\alpha}^{I}{}^{.i}\wedge
du_{I}^{\alpha}\wedge d^{n-1}x_{i}+d(\sum\nolimits_{|I|{}\leq k-1}p_{\alpha
}^{I}{}^{.i}u_{Ii}^{\alpha}-L[x,u])\wedge d^{n}x,
\]
where $d^{n-1}x_{i}:=i_{\partial/\partial x^{i}}d^{n}x$, and that the
corresponding PD-Hamilton equations (de Donder-Weyl equations) are nothing but
Euler-Lagrange equations determined by $\boldsymbol{S}_{1}$. We have thus
proved the following

\begin{theorem}
The de Donder-Weyl equations and the Euler-Lagrange equations determined by
$\boldsymbol{S}_{1}$ coincide.
\end{theorem}

Thus, despite the Legendre transform is far from being an isomorphism, the
Hamiltonian counterpart of the theory is basically identical to the Lagrangian one.

We conclude remarking that the geometric formulation of the (Hamilton-like,
higher derivative) de Donder field theory can be recovered from $\omega$
exactly as in \cite{v10}. This completes the program of the paper.

\end{document}